\documentclass[10pt,english]{smfart}

\usepackage[T1]{fontenc}
\usepackage[english]{babel}
\usepackage{amssymb,url,xspace,smfthm,MnSymbol,pdfsync,stmaryrd,mathtools,thmtools}
\usepackage{tikz,tikz-cd,denpastyle}

\newcommand{\BibTeX}{{\scshape Bib}\kern-.08em\TeX}
\newcommand{\T}{\S\kern .15em\relax }
\newcommand{\AMS}{$\mathcal{A}$\kern-.1667em\lower.5ex\hbox
{$\mathcal{M}$}\kern-.125em$\mathcal{S}$}

\DeclareMathOperator{\RHom}{R\calH om}

\DeclareMathOperator{\IC}{IC}
\DeclareMathOperator{\perH}{\prescript{\mathfrak{p}}{}{\mathcal{H}}}
\DeclareMathOperator{\perD}{\prescript{\mathfrak{p}}{}{D_{\parf}}}
\DeclareMathOperator{\BM}{BM}

\DeclareMathOperator{\parf}{ctf}

\DeclareMathOperator{\pH}{ {\prescript{\frakp}{} {\calH}}}
\DeclareMathOperator{\oH}{ {\prescript{o}{} {\calH}}}

\let\oldref\ref
\renewcommand{\ref}[1]{(\oldref{#1})}

\title{Motivic equivalence under semismall flops}
\author{\sc Wille Liu}

\address{Department of Mathematics \\
	National Taiwan University \\
Taipei}
\email{b01201052@ntu.edu.tw}

\begin{document}
\def\smfbyname{}

\begin{abstract}
	We prove that under semismall smooth flops, smooth projective varieties have (non-canonically) isomorphic Chow motives with coefficients in any noetherian local ring $\Lambda$ by comparing the pushforward of the constant intersection complexes through flopping contractions.
\end{abstract}
\begin{altabstract}
	Nous montrons que les vari\'et\'es projectives lisses ont les motifs de Chow \`a coefficient dans un anneau local n\oe th\'erien non-canoniquement isomorphes sous flops semi-petits par comparer les images directes des complexes d'intersection constants par contractions cr\'epantes.
\end{altabstract}
\maketitle

\tableofcontents

\setcounter{section}{-1}
\section{Introduction}
Certain strong relations between birational minimal models have long been suspected. Koll\'ar \cite{Kol} proved that 3-dimensional birational $\bfQ$-factorial terminal minimal models have isomorphic Hodge structures on intersection cohomology groups by employing intersection cohomology theory. For higher dimensions, it was proven by Batryev \cite{Bat} that birational Calabi-Yau manifolds have the same Betti numbers by using p-adic integration. This result was extended by Wang \cite{wang1} to the classes of K-equivalent smooth projective manifolds which in particular applies to birational smooth minimal models. With the aid of motivic integration developped by Denef and Loeser these results were further refined \cite{DL} to the equivalence of Hodge numbers. Most of the results nowadays are either established by certain integration formalism to obtain ``numerical results'' \cite{wang2}, or are restricted to particular cases of flops to obtain rather strong ``geometric conclusions'' \cite{LLW},\cite{FW}. A series of conjectures concerning K-equivalent proper smooth varieties are proposed in \cite[Section 6]{wang2} and explained in \cite[Section 4.3]{wang3}. The first one of these conjecture is to find an algebraic correspondence between K-equivalent projective smooth varieties that identifies their Chow motives. This article makes a progress in this direction.
\smallskip

The main result (\ref{main}) of this paper is of a style intermediating between these two ends, stating that two smooth projective varieties related by a \emph{semismall flop} (\autoref{sec2}) have isomorphic Chow motives with coefficients in a local ring. Following several ideas of \cite{Kol}, \cite{BM}, \cite{dCM1} and \cite{dCM2}, we adopt perverse sheaves to compare the motives of smooth varietes under semismall flops. 
\smallskip

For motives with coefficients $\bfQ$, the decomposition theorem is available and it suffices to compare local systems over respective strata (\ref{intcmp}). This can be achieved by employing arc spaces to compare irreducible components of the fibres. 
\smallskip

The decomposition theorem can fail, even for fields of positive characteristics. Nevertheless, this difficulty can be remedied for coefficients in local rings by the extensibility of morphisms between sheaves from an Zariski open subset (\ref{cat}) and then to the ability to lift invertibles in the category of finite associative algebras over a noetherian commutative local ring (\ref{lift}). These properties permit an abstract gluing of isomorphisms of perverse sheaves across strata. An immediate consequence (\ref{tor}) of the main result is that the singular cohomology groups $H^k(X, \bfZ)$ of a smooth variety $X$ are unchanged under semismall flops.
\smallskip

I am indebted to Professor Chin-Lung Wang for his supervision and innumerable inspiring discussions, and I am thankful to Chen-Yu Chi, Luc Illusie, Ming-Lun Hsieh, Hui-Wen Lin, Jeng-Daw Yu and many others for their advices.
\section{Preliminaries and notations}
\subsection{Derived categories of ctf complexes and perverse sheaves}
We will always work in the category of separated schemes of finite type over the complex numbers $\bfC$, an object whereof is simply called a variety henceforth. 
We will mostly compare complexes of sheaves in the bounded derived category $D^b_{\parf}(X, \Lambda) = D^b_{\parf}(X_{\mathrm{cl}}, \Lambda)$ of cohomologically constructible complexes sheaves of $\Lambda$ of finite tor-dimension on $X$ with coefficients in a (commutative unital) noetherian local ring $\Lambda$. See \cite{Sch} and \cite[expos\'e 2]{SGA4.5} for a reference. Henceforth we will simply call sheaves in place of sheaves of $\Lambda$-modules We summerise briefly here important properties of this category. 
\smallskip

Recall that an (algebraic) \emph{Whitney stratification} of an algebraic variety $X$ is a partition (called stratification) of $X$ into finite disjoint collection of Zariski locally closed subsets $\calT = \left\{ T_i \right\}_{i\in I}$, whereof each element (called stratum) $T$ is a smooth subvarieties such that the closure $\ba T\subseteq X$ is a union of strata. Every stratification admits a Whitney stratification as refinement.
\smallskip

According to \cite{BBD}, by a stractification of a algebraic variety $X$ we mean a Whitney stratification of equidimensional strata of $X$.
\smallskip

We call a sheaf $\calF$ of $\Lambda$-modules on  $X$ a \emph{locally constant constructible sheaf}, if there is an open covering  $\left\{ U_i \right\}_{i\in I}$ (in the classical topology) of $X$ such that the restriction of $\calF$ to each $U_i$ is a finitely generated free $A$-module. 
\smallskip

A sheaf $\calF$ is called \emph{constructible} if there is a Whitney stratification of $X$ along each stratum of which $\calF$ is locally constant constructible.  A complex of sheaves is called \emph{cohomologically constructible} if every cohomology sheaf of it is constructible. A complex of sheaves $K$ is called \emph{of finite torsion dimension} if each of its stalk $K_x$ is isomorphic in the bounded derived category of $\Lambda$-modules $D^b(\Lambda)$ to a complex of projective $\Lambda$-modules. 
\smallskip

Let $D^b_{\parf}(X, \Lambda)$ denote the full subcategory of $D^b(X, \Lambda)$ of bounded cohomologically constructible complexes of finite torsion dimension (ctf for short).
\smallskip

Amongst the most important topological properties of an algebraic map is the preservation of constructibility. Given an morphism of varieties $f\colon X\to Y$, the six operations are induced on the subcategories:
\begin{eq}
	f_*, f_!\colon D^b_{\parf}(X, \Lambda) \to D^b_{\parf}(Y, \Lambda),\quad f^*, f^!\colon D^b_{\parf}(Y, \Lambda) \to D^b_{\parf}(X, \Lambda),
\end{eq}
as well as the derived tensor product
\begin{eq}
	\relbar\otimes^L\relbar\colon D^b_{\parf}(X, \Lambda)\times D^b_{\parf}(X, \Lambda) \to D^b_{\parf}(X, \Lambda)
\end{eq}
and the internal homomorphism
\begin{eq}
	\RHom_X\left(\relbar, \relbar\right)\colon D^b_{\parf}(Y, \Lambda)^{\op}\times D^b_{\parf}(X, \Lambda) \to D^b_{\parf}(X, \Lambda).
\end{eq}
Proofs of the constructibility can be found in \cite[section 4.2.2]{Sch}.
As an immediate consequence, $R\Hom_X\left( K, L \right) \in D^b_{\parf}(\Spec \bfC, \Lambda)$ is a ctf $\Lambda$-complex for any pair $K, L\in D^b_{\parf}(X, \Lambda)$.
\smallskip

There is a natural biduality functor $\left( \relbar \right)^{\vee}\colon D^b_{\parf}(X, \Lambda)\to D^b_{\parf}(X, \Lambda)^{\op}$ defined by $K = \RHom\left( K, \omega_X  \right)$, where $\omega_X = \pi_X^!\Lambda$ is the dualising complex, and $\pi_X\colon X\to \Spec\bfC$ is the constant map. 

Perverse t-structures of middle perversity are defined on these categories of the following manner:
\begin{eq}
	\perD^{\le 0} = \left\{ K\in D^b_{\parf}(X, \Lambda) :  \left(\forall i\colon S\hookrightarrow X\; \text{stratum of $X$}\right)\left(\forall k > \dim S\right)\;\calH^k i^* K = 0  \right\}
\end{eq}
and
\begin{eq}
	\perD^{\ge 0} = \left\{ K\in D^b_{\parf}(X, \Lambda) :  \left(\forall i\colon S\hookrightarrow X\; \text{stratum of $X$}\right)\left(\forall k < \dim S\right)\;\calH^k i^! K = 0  \right\}.
\end{eq}
The heart $\calP(X, \Lambda)$ of this t-structure is simply called the category of perverse sheaves on $X$. We remark that $\calP(X, \Lambda)$ is a noetherian abelian category and that it is not in general artinian. 
\smallskip
\subsection{Arc spaces}
See \cite{DL} for detailed properties of arc spaces. To each variety $X$ we can associate a projective system of schemes $\scrL_m(X)$, called the \emph{$m$\textsuperscript{th} truncated arc space} of $X$, which represents the functor $\underline{\Hom}_{\Spec\bfC}\left( \Spec \bfC[t]/t^{m+1}, X \right)$ that sends a $\bfC$-scheme $T$ to $\Hom_{\Spec\bfC}\left(T\times_{\Spec \bfC} \Spec \bfC[t]/t^{m+1}, X \right)$. Let $\calL(X)$ be the limit of this system, which represents the functor $\underline{\Hom}_{\Spec\bfC}\left( \Spec \bfC\llbracket t\rrbracket, X \right)$.
When $X$ is smooth of equidimension $n$, the natural projection $\pi_{m, X}\colon \scrL_m(X) \to X$ is a Zariski locally trivial $\bfA^{mn}$-fibration over $X$. 
\smallskip

When there is a birational morphism $f\colon X\to Y$ between smooth varieties, one can read the structure of the system of induced morphisms $\scrL_m(f)\colon \scrL_m(X)\to \scrL_m(Y)$ from the relative canonical divisor $\calJ = K_{X/Y}$ and vice versa. 
We have an important lemma about the structure of a birational morphism~\cite[Lemma 3.4]{DL}.
\begin{lemm}\label{str}
	Let $X$ and $Y$ be smooth algebraic varieties over $k$, of pure dimension $d$ and let $f\colon X\to Y$ be a birational morphism. For $k$ in $\bfN$, let 
	\begin{eq}
		\scrL(X)^{k} \colon= \left\{ \gamma \in \scrL(X) : \left(\ord_t \mathcal{J}\right)(\gamma) = k \;\text{and}\; \scrL(f) \gamma \in \scrL(Y)\right\}.
	\end{eq}
	For $m\in \bfN$, let $\scrL_m(f)\colon \scrL_m(X) \to \scrL_m(Y)$ be the morphism induced by $f$, and let $\scrL_m(X)^k$ be the image of $\scrL(X)^k$ in $\scrL_m(X)$. Then, for all $k$ in $\bfN$ with $m\ge 2k$, the following holds.
	\begin{enumerate}
		\item[(a)]
			The set $\scrL_m(X)^k$ is a union of fibres of $\scrL_m(f)$.
		\item[(b)]
			The restriction of $\scrL_m(f)$ to $\scrL_m(X)^k$ is a piecewise trivial fibration with fibre $\bfA^k$ onto its image.
	\end{enumerate}
\end{lemm}
Besides $\scrL_m(X)^k$, we introduce the notations
\begin{eq}
	\scrL(X)^{\le k} \colon= \left\{ \gamma \in \scrL(X) : \left(\ord_t \mathcal{J}\right)(\gamma) \le k \;\text{and}\; \scrL(f) \gamma \in \scrL(Y)\right\}, \\
\scrL(X)^{> k} \colon= \left\{ \gamma \in \scrL(X) : \left(\ord_t \mathcal{J}\right)(\gamma) > k \;\text{and}\; \scrL(f) \gamma \in \scrL(Y)\right\},
\end{eq}
etc, to indicate the subvarieties of arcs of corresponding orders, and denote their images by $\scrL_m(X)^{\le k}$, $\scrL_m(Y)^{\le k} = \scrL_m(f)\left( \scrL_m(X)^{\le k} \right)$, etc.
\smallskip

Moreover, for a subset $S\subseteq X$, we add the restriction sign $\relbar\Big|_S$ to indicate those arcs $\gamma$ on $X$ that are originated from $S$, that is to say $\gamma(0)\in S$. For example, $\scrL_m(X)\Big|_S = \pi_{m, X}^{-1}(S) \subseteq \scrL_m(X)$. Notations concerning order of arcs in the previous paragraph also applies to this situation.

\section{Comparison over strata}\label{sec2}
In this section, we will consider complexes of sheaves of abelian groups and the bounded derived categories $D^b(X, \bfZ)$ thereof.
\smallskip


Let $X$ and $X'$ be smooth varieties of dimension $n$, $Y$ a varieties, and $f\colon X \to Y$ and $f'\colon X'\to Y$ be proper birational morphisms. Suppose that $X$ and $X'$ are K-equivalent through ${f'}^{-1} \circ f$. That is to say, given any common resolution $g\colon Z\to X$ and $g'\colon Z \to X'$, the relative canonical divisors $K_g$ and $K_{g'}$ are equal. We suppose further that $f$ and $f'$ are \emph{semismall} in the sense that $\dim X\times_Y X = n$ and $\dim X'\times_Y X' = n$. In this case ${f'}^{-1}\circ f$ is refered to as a \emph{semismall K-equivalence}.
\smallskip

One important example we shall bear in mind is that when $f\colon X\to Y$ is a semismall crepant resolution and $f'\colon X'\to Y$ is a flop, $X$ and $X'$ are K-equivalent through ${f'}^{-1} \circ f$. In this case ${f'}^{-1}\circ f$ is refered to as a \emph{semismall flop}.
\smallskip

A stratum $T\in \calT$ is called \emph{$f$-relevant} if $\dim f^{-1}(T)\times_{T} f^{-1}(T) = n$.
\smallskip

The goal of this section is to prove the following \ref{intcmp}.
\begin{prop}\label{intcmp}
	For a sufficiently fine stratification $\calT$ of $Y$, on each relevant stratum $T\in\calT$ of dimension $d$, there is an isomorphism $R^{n-d}f_*\bfZ_{f^{-1}(T)} \cong R^{n-d}f'_*\bfZ_{ {f'}^{-1}(T)}$. 
\end{prop}
Taking the dual local system,
\begin{coro}\label{dual}
	Under the assumptions of \ref{intcmp}, let $i\colon T\to Y$ be the embedding. There is an isomorphism
	$R^{n-d}i^!f_*\bfZ_{X} \cong R^{n-d}i^!f'_*\bfZ_{ X'}$
\end{coro}

Here is a plan of the proof. 
\smallskip

Given an initial stratification $\calT$, we shall recursively study the isomorphism and at the same time refine the stratification and the common resolution $Z$ if necessary. \\
To be more precise, assume the comparison is done on a union of strata, which is Zariski open in $Y$. We shall pick a stratum dense, say $T$, in the complement, blow up $Z$ accordingly, and establish the comparison on a Zariski open subset of $T$ in case $T$ is relevant, by means of arc spaces. Then we shall refine the stratification by splitting up $T$ into the Zariski open subset and some other smooth Zariski locally closed subvarieties. By noetherian induction, this will prove \ref{intcmp}.
\smallskip

Let $T\subseteq Y$ be a Zariski locally closed smooth conncected subvariety, of dimension $d$. Suppose that $R^{n-d} f_* \bfZ$ and $R^{n-d}f'_*\bfZ$ are (classically) locally constant constructible, and that $f$ and $f'$ are flat over $T$.
\smallskip

We begin with some observations on the irreducible components.  Consider the projection 
\begin{eq}
	\pi_{m, f^{-1}(T)}\colon \scrL_m(X)\Big|_{f^{-1}(T)} \to f^{-1}(T)
\end{eq}
from the space of $m$\textsuperscript{th} truncated arcs with origin in $f^{-1}(T)$, to $f^{-1}(T)$. 

\begin{lemm}\label{toarc}
	Under the above assumptions, the irreducible components of $\scrL_m(X)\Big|_{f^{-1}(T)}$ correspond to those of $f^{-1}(T)$ under $\pi_{m, f^{-1}(T)}$, and that
	\begin{eq}
		R\pi_{m, f^{-1}(T), !}\colon \bfZ[2mn] \cong \bfZ.
	\end{eq}
\end{lemm}
\begin{proof}
	Let $E\subset f^{-1}(T)$ be an irreducible component. Then $\scrL_m(X) \Big|_E$ is a Zariski locally trivial $\bfA^{mn}$-fibration over $E$, and hence is irreducible. Consequently, the irreducible components of $\scrL_m(X)\Big|_{f^{-1}(T)}$ correspond to those of $f^{-1}(T)$ under the projection $\pi_{f^{-1}(T), m}\colon \scrL_m(X)\Big|_{f^{-1}(T)} \to f^{-1}(T)$. Since it is a $\bfA^{mn}$-fibration, by the base change property for exceptional push-forwards, the projection $\pi_{m, f^{-1}(T)}$ induces
	\begin{eq}
		R\pi_{m, f^{-1}(T), !}\colon \bfZ[2mn] \cong \bfZ.
	\end{eq}
\end{proof}
In regards to \ref{toarc}, we have
\begin{eq}
	R\pi_{m, f^{-1}(T), !}\colon \bfZ[2mn] \cong \bfZ \quad \text{on $f^{-1}(T)\subseteq X$}
\end{eq}
as well as 
\begin{eq}
	R\pi_{m, {f'}^{-1}(T), !}\colon \bfZ[2mn] \cong \bfZ\quad \text{on ${f'}^{-1}(T)\subseteq X'$}
\end{eq}
and 
\begin{eq}
	R\pi_{m, h^{-1}(T), !}\colon \bfZ[2mn] \cong \bfZ\quad \text{on $h^{-1}(T)\subseteq Z$}.
\end{eq}
The situation is indicated below
\begin{eq}
	\begin{tikzcd}
		& \scrL_m(Z)\Big|_{h^{-1}(T)} \arrow{dl}{\scrL_m(g)}\arrow{d}{\pi_{m, h^{-1}(T)}}\arrow{dr}{\scrL_m(g')} & \\
		\scrL_m(X)\Big|_{f^{-1}(T)}\arrow{d}{\pi_{m, f^{-1}(T)}} & h^{-1}(T)\arrow{dl}{g}\arrow{dr}{g'} & \scrL_m(X')\Big|_{{f'}^{-1}(T)}\arrow{d}{\pi_{m, {f'}^{-1}(T)}} \\
		f^{-1}(T) &  & {f'}^{-1}(T) \\
	\end{tikzcd}
\end{eq}
\smallskip

The next step is to establish relations between the irreducible components of $\scrL_m(Z)\Big|_{h^{-1}(T)}$ and those of $\scrL_m(X)\Big|_{f^{-1}(T)}$. We shall see that each component of $\scrL_m(X)\Big|_{f^{-1}(T)}$ is dominated by one certain irreducible component of $\scrL_m(Z)\Big|_{h^{-1}(T)}$. 
\smallskip

For each irreducible component $F \subseteq \scrL_m(Z)\Big|_{h^{-1}(T)}$, let 
\begin{eq}
	\delta(F) = \min\left\{k : \scrL_m(Z)^{\le k}\cap F \neq \emptyset  \right\}.
\end{eq}
\begin{lemm}\label{cmpcomp}
Under the above assumptions, assume in addition that $m \ge 2\delta(F) + 1$ for all irreducible components $F\subseteq \scrL_m(Z)\Big|_{h^{-1}(T)}$. Then $\delta(F)  = \dim F - \dim \scrL_m(g)(F) \ge \dim F - (n+d)/2 + mn$, and that the equality holds exactly when $F$ dominates an irreducible component of maximal dimension $(n + d)/2 + mn$.
\end{lemm}
\begin{proof}
	\smallskip

	By \ref{str}, the additional assumption $m \ge 2\delta(F) + 1$ make arcs $\scrL(Z)^{\le \delta(F)}$ of contact order $\le \delta(F)$ all stabilised(i.e. there is a piecewise affine fibration structure), in particular, $\scrL_m(g)$ is a Zariski locally trivial $\bfA^{\delta(F)}$-fibration when restricted to $\scrL_m(Z)^{\le \delta(F)}$ with image $\scrL_m(X)^{\le \delta(F)}$. Therefore 
	\begin{eq}
		\dim\left(\scrL_m(Z)^{\le \delta(F)} \cap F\right)- \dim\left(\scrL_m(X)^{\le \delta(F)} \cap \scrL_m(g)(F)\right) = \delta(F).
	\end{eq}
	\smallskip

	On the other hand, $\scrL_m(Z)^{\le \delta(F)}\cap F$ is Zariski open and dense in $F$, and that $\scrL_m(X)^{\le \delta(F)}\cap \scrL_m(g)(F)$ is Zariski open and dense in $\scrL_m(g)(F)$, so 
	\begin{eq}
		\delta(F) &= \dim\left(\scrL_m(Z)^{\le \delta(F)} \cap F\right) - \dim\left(\scrL_m(X)^{\le \delta(F)} \cap \scrL_m(g)(F)\right) \\
		&= \dim F - \dim \scrL_m(g)(F).
	\end{eq}
	
	Recall that the assumption on semismallness forces that $\dim f^{-1}(T) \le (n + d)/2$, and hence $\dim \scrL_m(g)(F)\le \dim \scrL_m(X)\Big|_{f^{-1}(T)}\le (n + d)/2 + mn$. 
	\begin{eq}
		\delta(F) = \dim F - \dim \scrL_m(g)(F) \ge \dim F - \frac{n + d}2 - mn.
	\end{eq}
	The equality holds exactly when $\dim \scrL_m(g)(F) = (n + d)/2 - mn$. In that case, $\scrL_m(g)(F)$ is of maximal dimension in $\scrL_m(X)\Big|_{f^{-1}(T)}$.
\end{proof}
\begin{coro}\label{cmpcomp1}
	Under the assumptions of \ref{cmpcomp}, $\scrL_m(g)$ induces a bijection
	\begin{eq}
		&\left\{ \text{irreducible components $E\subseteq \scrL_m(X)\Big|_{f^{-1}(T)}$} : \dim E = \frac{n+d}2 + mn \right\} \\
		\longleftrightarrow &\left\{ \text{irreducible components $F\subseteq \scrL_m(Z)\Big|_{h^{-1}(T)}$} : \delta(F) = \dim F - \frac{n + d}2 - mn \right\}
	\end{eq}
\end{coro}
\begin{proof}
	In regards to \ref{cmpcomp}, we are left to prove that each irreducible component $E\subseteq \scrL_m(X)\Big|_{f^{-1}(T)}$ is dominated by exactly one irreducible component of $\scrL_m(Z)\Big|_{h^{-1}(T)}$. This is clear by the surjectivity of $\scrL_m(g)$ and the fact that $\scrL_m(g)$ is an $\bfA^{\delta(F)}$-fibration over an Zariski open set of $E$. 
\end{proof}
For convenience, we may and we shall refine the resolution as $Z' \to Z$, which is an isomorphism outside the closure $\ba{h^{-1}(T)}$, such that $h^{-1}(T)$ is a sum of divisors, and that $\dim F = (n - 1) + mn$ for all $F$. The condition then becomes $\delta(F) \ge (n - d)/2 - 1 =: \delta$. We rename $Z := Z'$. 
\smallskip

\ref{cmpcomp1} describes explicitly the relation between components of importance, the difficulty to the comparison being that $\scrL_m(g)$ is only an $\bfA^{\delta}$-fibration on an Zariski open subset. This difficulty is evitable by restricting everything to an Zariski open subset of $T$. 
\smallskip

Set
\small
\begin{eq}
	U = \left\{ t\in T : \dim\left(\scrL_m(X)^{> \delta}\Big|_{f^{-1}(t)}\right) < \frac{n-d}2 + mn,\quad  \dim\left(\scrL_m(X')^{> \delta}\Big|_{ {f'}^{-1}(t)}\right) < \frac{n-d}2 + mn\right\}.
\end{eq}
\normalsize

\begin{lemm}\label{res}
	Under the assumptions of \ref{cmpcomp}, $U$ is Zariski open and dense in $T$.
\end{lemm}
\begin{proof}
	Put
	\begin{eq}
		V &= \left\{ t\in T : \dim\left(\scrL_m(X)^{> \delta}\Big|_{f^{-1}(t)}\right) < \frac{n-d}2 + mn \right\} \;\text{and}\\ 
		V' &= \left\{ t\in T : \dim\left(\scrL_m(X')^{> \delta}\Big|_{ {f'}^{-1}(t)}\right) < \frac{n-d}2 + mn \right\},
	\end{eq}
	so that $U = V \cap V'$. Then it suffices to verify that $V$ and $V'$ are Zariski open and dense in $T$. We shall do that for $V$ only.
	\smallskip

	We may assume that $T$ is irreducible.
	\smallskip

	Firstly, $V$ is nonempty. If $V = \emptyset$, then by the piecewise fibration structure \ref{str} $\dim\left( \scrL_m(Z)^{> \delta}\Big|_{h^{-1}(t)} \right) \ge \delta + 1 + (n-d)/2 + mn$ for all $t$, so 
	\begin{eq}
		\dim\left( \scrL_m(Z)^{> \delta}\Big|_{h^{-1}(T)} \right) \ge \delta + 1 + \frac{n-d}2 + mn + d = (m+1)n,
	\end{eq}
	which is absurd. Thus $V \neq \emptyset$.
	\smallskip

	We study the set $V$ by applying $R^{(n-d)+2mn}\left( f\circ \pi_{m, f^{-1}(T)} \right)_!$ to the sequence
	\begin{eq}
		0 \to \bfZ_{\scrL_m(X)^{\le \delta}\Big|_{f^{-1}(T)}} \to \bfZ_{\scrL_m(X)\Big|_{f^{-1}(T)}}\to \bfZ_{\scrL_m(X)^{> \delta}\Big|_{f^{-1}(T)}} \to 0.
	\end{eq}
	Then
	\small
	\begin{eq}
		R^{(n-d)+2mn}\left( f\circ \pi_{m, f^{-1}(T)} \right)_!\bfZ_{\scrL_m(X)\Big|_{f^{-1}(T)}}\to R^{(n-d)+2mn}\left( f\circ \pi_{m, f^{-1}(T)} \right)_!\bfZ_{\scrL_m(X)^{> \delta}\Big|_{f^{-1}(T)}} \to 0.
	\end{eq}
	\normalsize
	By \ref{toarc}, the former sheave is
	\begin{eq}
		R^{(n-d)+2mn}\left( f\circ \pi_{m, f^{-1}(T)} \right)_!\bfZ_{\scrL_m(X)\Big|_{f^{-1}(T)}}\cong R^{n-d}f_*\bfZ_{f^{-1}(T)},
	\end{eq}
	which is locally constant constructible by assumption. According to the ``semicontinuity'' of quotients of locally constant sheaves, the set of $t\in T$ at which the stalk is zero
	\begin{eq}
		0 = \left(R^{(n-d)+2mn}\left( f\circ \pi_{m, f^{-1}(T)} \right)_!\bfZ_{\scrL_m(X)^{> \delta}\Big|_{f^{-1}(T)}}\right)_t = H^{(n-d) + 2mn}_c\left( \scrL_m(X)^{> \delta}\Big|_{f^{-1}(t)}, \bfZ \right)
	\end{eq}
	is Zariski open in $T$. Since $\scrL_m(X)\Big|_{f^{-1}(T)} \le d + mn$, $H^{(n-d) + 2mn}_c\left( \scrL_m(X)^{> \delta}\Big|_{f^{-1}(T)}, \bfZ \right) = 0$ exactly when $\dim\left( \scrL_m(X)^{> \delta}\Big|_{f^{-1}(T)}\right) < d + mn$. The set is exactly $V$, whence $V$ is Zariski open.
\end{proof}

By restricting to $\scrL_m(Z)^{\le \delta}\Big|_{h^{-1}(U)}$ and its image $U\subseteq T$, we arrive at the new situation indicated below:  
\begin{eq}
	\begin{tikzcd}
		& \scrL_m(Z)^{\le \delta}\Big|_{h^{-1}(U)} \arrow{dl}{\scrL_m(g)^{\le \delta}}\arrow{d}{\pi_{m, h^{-1}(U)}}\arrow{dr}{\scrL_m(g')^{\le \delta}} & \\
		\scrL_m(X)^{\le \delta}\Big|_{f^{-1}(U)}\arrow{d}{\pi_{m, f^{-1}(U)}} & h^{-1}(U)\arrow{dl}{g}\arrow{dd}{h}\arrow{dr}{g'} & \scrL_m(X')^{\le \delta}\Big|_{{f'}^{-1}(U)}\arrow{d}{\pi_{m, {f'}^{-1}(U)}} \\
		f^{-1}(U) \arrow{dr}{f} &  & {f'}^{-1}(U)\arrow{dl}{f'} \\
		& U &  \\
	\end{tikzcd}
\end{eq}
According to the structure theorem (\ref{str}), the morphisms $\scrL_m(g)^{\le \delta}\colon \scrL_m(Z)^{\le \delta}\Big|_{h^{-1}(U)}\to \scrL_m(X)^{\le \delta}\Big|_{f^{-1}(U)}$ and $\scrL_m(g')^{\le \delta} \colon \scrL_m(Z)^{\le \delta}\Big|_{h^{-1}(U)} \to \scrL_m(X')^{\le \delta}\Big|_{{f'}^{-1}(U)}$ are piecewise trivial $\bfA^{\delta}$-fibrations. Therefore there are Zariski open dense $V\subseteq \scrL_m(X)^{\le \delta}\Big|_{f^{-1}(U)}$ and $V'\subseteq \scrL_m(X)^{\le \delta}\Big|_{{f'}^{-1}(U)}$ such that $\scrL_m(g)^{\le \delta}$ and $\scrL_m(g')^{\le \delta}$ are trivial $\bfA^{\delta}$-fibrations. 
Let $V^c\subseteq \scrL_m(X)^{\le \delta}\Big|_{f^{-1}(U)}$ and ${V'}^c\subseteq \scrL_m(X')^{\le \delta}\Big|_{{f'}^{-1}(U)}$ be the corresponding complements.
\smallskip



Finally, pick $W\subseteq U$ such that the fibre of $f\circ \pi_{m, f^{-1}(U)}\colon V^c \to U$ and of $f'\circ \pi_{m, {f'}^{-1}(U)}\colon {V'}^c \to U$ over each point of $W$ is of dimension $< (n-d) / 2 + mn$.

\begin{lemm}\label{cohomcmp}
	Under the assumptions of \ref{cmpcomp}, in the above diagramme, there is an isomorphism on $W$, $R^{n-d}f_*\bfZ_{f^{-1}(W)} \cong R^{n - d}f'_*\bfZ_{{f'}^{-1}(W)}$.
\end{lemm}
\begin{proof}
	In this context, Leray spectral sequence for $R\pi_{m, f^{-1}(U), !}\circ R\scrL_m(g)^{\le \delta}_!$ then gives
	\begin{align}\label{eq1}
		R^{2mn + 2\delta}\left(\pi_{m, f^{-1}(U)} \circ \scrL_m(g)^{\le \delta}\right)_! \bfZ_{\scrL_m(Z)^{\le\delta}} \cong  R^{2mn}\pi_{m, f^{-1}(U), !}\left(R^{2\delta}\scrL_m(g)^{\le \delta}_! \bfZ_{\scrL_m(Z)^{\le\delta}}\right) 
	\end{align}
	\smallskip
	
	Since $\scrL_m(g)$ is trivial $\bfA^{\delta}$-fibration over $V$. 
	\begin{align}\label{eq2}
		R^{2\delta}\scrL_m(g)^{\le \delta}_! \bfZ_{\scrL_m(Z)^{\le\delta}}\Big|_V \cong \bfZ_V
	\end{align}
	\smallskip

	Applying $R\left( f\circ \pi_{m, f^{-1}(U)}\right)_!$ to the relative sequence on $\scrL_{m}(X)\Big|_{f^{-1}(U)}$
	\begin{eq}
		0\to \bfZ_{V} \to \bfZ_{\scrL_m(X)^{\le \delta}\Big|_{f^{-1}(U)}} \to \bfZ_{V^c} \to 0, 
	\end{eq}
	we have
	\begin{align}\label{eq3}
		\left(R^{(n-d) + 2mn}\left( f\circ \pi_{m, f^{-1}(U)}\right)_! \bfZ_{V}\right)\Big|_{W} \cong \left(R^{(n-d) + 2mn}\left( f\circ \pi_{m, f^{-1}(U)}\right)_! \bfZ_{\scrL_m(X)^{\le \delta}\Big|_{f^{-1}(U)}}\right)\Big|_{W}.
	\end{align}
	\smallskip

	Applying $R\left( f\circ \pi_{m, f^{-1}(U)}\right)_!$ to the relative sequence on $\scrL_{m}(X)\Big|_{f^{-1}(U)}$
	\begin{eq}
		0\to \bfZ_{\scrL_m(X)^{\le \delta}\Big|_{f^{-1}(U)}} \to \bfZ_{\scrL_m(X)\Big|_{f^{-1}(U)}} \to \bfZ_{\scrL_m(X)^{> \delta}\Big|_{f^{-1}(U)}} \to 0, 
	\end{eq}
	we have then
	\begin{align}\label{eq4}
		R^{(n-d) + 2mn}\left( f\circ \pi_{m, f^{-1}(U)}\right)_! \bfZ_{\scrL_m(X)^{\le \delta}\Big|_{f^{-1}(U)}} \cong R^{(n - d) + 2mn}\left( f\circ \pi_{m, f^{-1}(U)}\right)_! \bfZ_{\scrL_m(X)\Big|_{f^{-1}(U)}}.
	\end{align}
	and the latter is isomorphic to
	\begin{align}\label{eq5}
		R^{n - d}f_*\left( R^{2mn}\pi_{m, f^{-1}(U), !} \bfZ_{\scrL_m(X)\Big|_{f^{-1}(U)}}\right) \cong R^{n - d}f_* \bfZ_{f^{-1}(U)}.
	\end{align}
	\smallskip

	Combining equations \ref{eq1}, \ref{eq2}, \ref{eq3}, \ref{eq4} and \ref{eq5}, we get
	\begin{eq}
		&\left(R^{(n - d) + 2mn+ 2\delta}\left( h\circ \pi_{m, h^{-1}(U)} \right)_!\bfZ_{\scrL_m(Z)^{\le \delta}\Big|_{h^{-1}(U)}}\right)\Big|_{W} \\
		&\cong \left(R^{(n - d) + 2mn+ 2\delta}\left( f\circ\pi_{m, f^{-1}(U)}\circ \scrL_m(g) \right)_!\bfZ_{\scrL_m(Z)^{\le \delta}\Big|_{h^{-1}(U)}}\right)\Big|_{W} \\
		&\cong \left(R^{(n - d) + 2mn}\left( f\circ\pi_{m, f^{-1}(U)}\right)_!\bfZ_{V}\right)\Big|_{ W} \\
		&\cong \left(R^{(n - d) + 2mn}\left( f\circ\pi_{m, f^{-1}(U)}\right)_!\bfZ_{\scrL_m(X)^{\le \delta}\Big|_{f^{-1}(U)}}\right)\Big|_{ W} \\
		&\cong \left(R^{(n - d) + 2mn}\left( f\circ\pi_{m, f^{-1}(U)}\right)_!\bfZ_{\scrL_m(X)\Big|_{f^{-1}(U)}}\right)\Big|_{ W} \\
		&\cong \left( R^{n - d}f_* \bfZ_{f^{-1}(U)}\right)\Big|_{ W} \\
	\end{eq}
	Similarly, 
	\begin{eq}
		\left(R^{(n - d) + 2mn+ 2\delta}\left( h\circ \pi_{m, h^{-1}(U)} \right)_!\bfZ_{\scrL_m(Z)^{\le \delta}\Big|_{h^{-1}(U)}}\right)\Big|_{ W} \cong\left( R^{n - d}f'_* \bfZ_{{f'}^{-1}(U)}\right)\Big|_{ W}.
	\end{eq}
	Therefore 
	\begin{eq}
		R^{n - d}f_* \bfZ_{f^{-1}(W)}\cong \left( R^{n - d}f_* \bfZ_{f^{-1}(U)}\right)\Big|_{W} \cong \left( R^{n - d}f'_* \bfZ_{{f'}^{-1}(U)}\right)\Big|_{W} \cong R^{n - d}f'_* \bfZ_{{f'}^{-1}(W)}.
	\end{eq}
	We may redefine $T$ as $W$ and refine the stratification $\calT$ accordingly, so that
	\begin{eq}
		R^{n - d}f_* \bfZ \cong R^{n - d}{f'}_* \bfZ.
	\end{eq}
\end{proof}
From the arguments above, one sees that 
\begin{prop}\label{dim}
	Under the assumptions of \ref{cmpcomp}, the semismallness of $f$ implies that of $f'$, and that $\dim X\times_Y X' = n$.
	Moreover, if a stratum $T$ is $f$-relevant, then it is also $f'$-relevant.
\end{prop}
Finally, we can prove \ref{intcmp}. 
\begin{proof}[Proof of \ref{intcmp}]
	As is sketched right after the statement of \ref{intcmp}, the proof proceeds by recurrence on strata. 
	\smallskip

	We assign a partial order on the set of strata. There is a relation between two strata $T \le T'$ if $T \subseteq \ba{T'}$.
	\smallskip

	Suppose that $T\in \calT$ is a stratum, and that on every $T'\in \calT$ with $T \leq T'$ but $T\neq T'$,  the comparison $R^{n - d}f_*\bfZ_T \cong R^{n - d}{f'}_*\bfZ_T$ has been established. Apply \ref{cohomcmp} to $T$, so that there is a resolution $Z_1 \to Z$ and an Zariski open $\ba W\subseteq T$ with $R^{n - d}f_*\bfZ_{\ba W} \cong R^{n - d}{f'}_*\bfZ_{\ba W}$. Refine the stratification on $\ba T$ so that $\ba W$ is one of the strata. This finishes the recursive step.
	\smallskip

	Clearly by the noetherian assumption on $Y$, the recurrence eventually stops. That finishes the proof.
\end{proof}
\begin{prop}\label{Qmot}
	In coefficient $\bfQ$, we have
	\begin{eq}
		Rf_* \bfQ \cong Rf'_* \bfQ
	\end{eq}
	Consequently, if $X$ and $X'$ are projective, then the Chow motives in $\bfQ$-coefficients are isomorphic
	\begin{eq}
		M(X)_{\bfQ} \cong M(X')_{\bfQ}
	\end{eq}
\end{prop}
\begin{proof}
	Tensoring with $\bfQ$ the isomorphisms provided by \ref{intcmp}, we have
	\begin{eq}
		R^{n - d}f_*\bfQ &\cong R^{n - d}f_*\left(\bfZ\otimes^L f^*\bfQ\right)\cong R^{n - d}f_*\bfZ\otimes^L \bfQ \\
		&\cong R^{n - d}f'_*\bfZ\otimes^L \bfQ \cong R^{n - d}f'_*\left(\bfZ\otimes^L {f'}^*\bfQ\right) \cong R^{n - d}f'_*\bfQ.
	\end{eq}
	Then 
	\begin{eq}
		\IC\left(T, R^{n-d}f_*\bfQ \right) \cong \IC\left(T, R^{n-d}f'_*\bfQ \right).
	\end{eq}
	By the explicit decomposition theorem of \cite{BM}, 
	\begin{eq}
		Rf_*\bfQ[n] \cong \bigoplus_{T\in \calT}\IC\left(T, R^{n-d}f_*\bfQ \right) \cong \bigoplus_{T\in \calT}\IC\left(T, R^{n-d}f'_*\bfQ \right) \cong Rf'_*\bfQ[n] .
	\end{eq}
	An argument in \cite{dCM2}, which we will reproduce in the next section, shows that $M(X)_{\bfQ}\cong M(X')_{\bfQ}$.
\end{proof}
As the decomposition theorem is available only with coefficients in fields of characteristic $0$, in order to extend the result to more general coefficients, we shall study the extensions of perverse sheaves across strata in the next section and forth.

\section{Recollement of perverse sheaves}
Provided with \ref{intcmp} we are left with extending the isomorphisms across strata, in other words, to glue the isomorphisms of perverse sheaves over strata. The gluing is not automatic. For nice general expositions of gluing of perverse sheaves, confer \cite{BBD} or \cite{Jut}. 
\smallskip

Now fix a noetherian coefficient ring $\Lambda$. Throughout this section, we work on the derived categories $D^b_{\parf}(Y, \Lambda)$ of ctf $\Lambda$-complexes with respect to a stratification $\calT$ of $Y$ obtained in \ref{intcmp}. According to the notations of \cite{BBD}, we will write $f_*, f^*, f_!, f^!$ for the derived functor, whereas standard cohomology sheaf functors will be denoted by $\oH^k$, and the perverse cohomology sheafs by $\perH^k$. 
\smallskip

The situation of semismall K-equivalence $f'^{-1}\circ f\colon X\dashrightarrow X'$ is continued in this section. According to \cite{BM}, the sheaves $f_*\Lambda[n]$ and $f'_*\Lambda[n]$ are perverse. That can be shown by estimation of supports using Leray spectral sequences together with the self-duality of $f_*\Lambda[n]$ and of $f'_*\Lambda[n]$.
\smallskip

We begin with a observation on the convolution algebra and morphisms between $f_*\Lambda[n]$ and $f'_*\Lambda[n]$.
\begin{lemm}[{\cite[Lemma 8.6.1]{CG}}]
	For $U\subseteq Y$ open subset, 
	\begin{eq}
		\Hom_{D^b_c(U, \Lambda)}\left( f_*\Lambda_{f^{-1}(U)}, {f'}_*\Lambda_{f^{-1}(U)} \right)\cong H^{\BM}_{2n}\left(f^{-1}(U) \times_{U} {f'}^{-1}(U), \Lambda\right). 
	\end{eq}
	and this isomorphism commutes with the obvious restriction morphisms of open subsets $V \subseteq U$ in the sense that the following diagramme commutes
	\begin{eq}
		\begin{tikzcd}
			\Hom_{D^b_c(U, \Lambda)}\left( f_*\Lambda_{f^{-1}(U)}, {f'}_*\Lambda_{{f'}^{-1}(U)} \right)\arrow{r}{\cong}\arrow{d}{\res} & H^{\BM}_{2n}\left(f^{-1}(U) \times_{U} {f'}^{-1}(U), \Lambda\right)\arrow{d}{\res} \\
			\Hom_{D^b_c(V, \Lambda)}\left( f_*\Lambda_{f^{-1}(V)}, {f'}_*\Lambda_{{f'}^{-1}(V)} \right)\arrow{r}{\cong} & H^{\BM}_{2n}\left(f^{-1}(V) \times_{V} {f'}^{-1}(V), \Lambda\right) 
		\end{tikzcd}
	\end{eq}
	When there is a third $f''\colon X'' \to Y$, the composition pairing
	\begin{eq}
		\begin{tikzcd}
			\Hom_{D^b_c(U, \Lambda)}\left( f_*\Lambda_{f^{-1}(U)}, {f'}_*\Lambda_{ {f'}^{-1}(U)} \right)\otimes_\Lambda\Hom_{D^b_c(U, \Lambda)}\left( f'_*\Lambda_{ {f'}^{-1}(U)}, f''_*\Lambda_{ {f''}^{-1}(U)} \right)\arrow{d}{\circ} \\
			\Hom_{D^b_c(U, \Lambda)}\left( f_*\Lambda_{f^{-1}(U)}, f''_*\Lambda_{ {f''}^{-1}(U)}\right)
		\end{tikzcd}
	\end{eq}
	is compatible with the convolution product of cohomological correspondences
	\begin{eq}
		\begin{tikzcd}
			H^{\BM}_{2n}\left(f^{-1}(U) \times_{U} {f'}^{-1}(U), \Lambda\right)\otimes_\Lambda H^{\BM}_{2n}\left( {f'}^{-1}(U) \times_{U} {f''}^{-1}(U), \Lambda\right)\arrow{d}{\bullet} \\
			H^{\BM}_{2n}\left( f^{-1}(U) \times_{U} {f''}^{-1}(U), \Lambda\right)
		\end{tikzcd}
	\end{eq}
\end{lemm}
Emphasis is put on the compatibility with restriction maps, whereof we will be using in a crucial way. Observe also that in the situation of semismall K-equivalence, since the Borel-Moore homology groups are freely generated by the respective sets of irreducible components of dimension $n$ (\ref{dim}), the restriction map is surjective, and so is the other restriction map. Besides, the compatibility of product structure provides a link between isomorphism of sheaves and of Chow motives.
\smallskip

\begin{coro}[\cite{dCM2}] \label{corr}
	In the situation of semismall K-equivalence ${f'}^{-1}\circ f\colon X\dashrightarrow X'$, any isomorphism $f_*\Lambda[n]\cong f'_*\Lambda[n]$ gives rise to a cohomological correspondence $\Gamma\in H^{\BM}_{2n}(X\times_Y X', \Lambda)$, which is the class of an algebraic cycle with coefficients in $\Lambda$. When $X$ and $X'$ are projective, $\Gamma$ provides an isomorphism of Chow motives $M(X)_\Lambda \cong M(X')_\Lambda$ with coefficients in $\Lambda$.
\end{coro}

We shall prove by recurrence on strata in $\calT$ that $f_*\Lambda[n] \cong f'_*\Lambda[n]$. We are then reduced to the situation under which there are Zariski open subsets $j\colon V\hookrightarrow U \subseteq Y$ with smooth complement $i\colon S = U\setminus V\hookrightarrow U$ of equidimension $d$, on which $\oH^{-d}f_* \Lambda[n]$ and  $\oH^{-d}f'_* \Lambda[n]$ are locally constant constructible sheaves, isomorphisms $\varphi\colon f_*\Lambda_{f^{-1}(V)}[n] \cong f'_*\Lambda_{ {f'}^{-1}(V)}[n]$ and $\oH^{-d}f_* \Lambda[n] \cong\oH^{-d}f'_* \Lambda[n]$.
\smallskip

Lift the isomorphism $\varphi$ to a morphism $\til\varphi\colon f_*\Lambda_{f^{-1}(U)}[n] \to f'_*\Lambda_{ {f'}^{-1}(U)}[n]$, and its inverse $\psi = \varphi^{-1}$ to $\til\psi\colon  f'_*\Lambda_{ {f'}^{-1}(U)}[n] \to f_*\Lambda_{f^{-1}(U)}[n]$. We remark that $\til \varphi$ and $\til \psi$ may not be isomorphisms anymore. 
\smallskip

Then there are morphisms between distinguished triangles
\begin{eq}
	\begin{tikzcd}
		i_*i^! f_* \Lambda[n] \arrow{r}\arrow{d} & f_*\Lambda[n] \arrow{r}\arrow{d}{\til\varphi} & j_*j^* f_*\Lambda[n] \arrow{r}\arrow{d}{\varphi} & i_*i^! f_* \Lambda[n+1]\arrow{d} \\
		i_*i^! f'_* \Lambda[n] \arrow{r}\arrow{d} & f'_*\Lambda[n] \arrow{r}\arrow{d}{\til\psi} & j_*j^* f'_*\Lambda[n] \arrow{r}\arrow{d}{\psi} & i_*i^! f'_* \Lambda[n+1]\arrow{d}\\
		i_*i^! f_* \Lambda[n] \arrow{r} & f_*\Lambda[n] \arrow{r} & j_*j^* f_*\Lambda[n] \arrow{r} & i_*i^! f_* \Lambda[n+1] \\
	\end{tikzcd}
\end{eq}
Recall that by the semismallness of $f$ and of $f'$, $f_*\Lambda[n]$ and $f'_*\Lambda[n]$ are perverse. Applying the perverse cohomolology functor $\pH^0$ to this diagramme, we have exact sequences
\begin{eq}
	\begin{tikzcd}[column sep=small]
		0\arrow{r} & \pH^0 i_*i^! f_* \Lambda[n] \arrow{r}\arrow{d} & f_*\Lambda[n] \arrow{r}\arrow{d}{\til\varphi} & \pH^0 j_*j^* f_*\Lambda[n] \arrow{r}{d}\arrow{d}{\varphi} & \pH^1 i_*i^! f_* \Lambda[n] \arrow{r}\arrow{d}{\ba \varphi} & 0\\
		0\arrow{r} & \pH^0 i_*i^! f'_* \Lambda[n] \arrow{r}\arrow{d} & f'_*\Lambda[n] \arrow{r}\arrow{d}{\til\psi} & \pH^0 j_*j^* f'_*\Lambda[n] \arrow{r}{d'}\arrow{d}{\psi} & \pH^1 i_*i^! f'_* \Lambda[n]\arrow{r}\arrow{d}{\ba \psi} & 0\\
		0\arrow{r} & \pH^0 i_*i^! f_* \Lambda[n] \arrow{r} & f_*\Lambda[n] \arrow{r} & \pH^0 j_*j^* f_*\Lambda[n] \arrow{r}{d} & \pH^1 i_*i^! f_* \Lambda[n]\arrow{r} & 0 \\
	\end{tikzcd}
\end{eq}
Since $\psi\varphi = 1$, $\ba\psi\ba\varphi = 1$ as well. Similarly $\ba\varphi\ba\psi = 1$, so $\ba\varphi$ and $\ba\psi$ are inverse to each other. The diagramme is reduced to 
\begin{eq}
	\begin{tikzcd}
		0\arrow{r} & \pH^0 i_*i^! f_* \Lambda[n] \arrow{r}\arrow{d} & f_*\Lambda[n] \arrow{r}\arrow{d}{\til\varphi} & \ker d \arrow{d}{\varphi}\arrow{r} & 0 \\
		0\arrow{r} & \pH^0 i_*i^! f'_* \Lambda[n] \arrow{r}\arrow{d} & f'_*\Lambda[n] \arrow{r}\arrow{d}{\til\psi} & \ker d' \arrow{d}{\psi}\arrow{r} & 0 \\
		0\arrow{r} & \pH^0 i_*i^! f_* \Lambda[n] \arrow{r} & f_*\Lambda[n] \arrow{r} & \ker d \arrow{r} & 0  \\
	\end{tikzcd}
\end{eq}
where $\varphi$ and $\psi$ are isomorphisms. Since $\pH^0 i_*i^! f_* \Lambda[n] \cong \left(\oH^{n-d} i_*i^! f_* \Lambda\right)[d] \cong \left(\oH^{n-d} i_*i^! f'_* \Lambda\right)[d] \cong \pH^0 i_*i^! f'_* \Lambda[n]$ by \ref{dual} tensored with $\Lambda$, it amounts now to a statement concerning extensions of perverse sheaves, which will be demonstrated in the next section. 
\begin{lemm}\label{cat}
	Let $A$ and $C$ be objects in an abelian category $\calC$ locally finite over a noetherian local ring $\Lambda$, such that $\Hom_{\calC}(A, C) = 0$. Suppose we have a diagramme of short exact sequences
	\begin{eq}
		\begin{tikzcd}
			0\arrow{r} & A \arrow{r}{\alpha}\arrow{d} & B \arrow{r}{\beta}\arrow{d}{\til\varphi} & C\arrow{d}{\varphi}\arrow{r} & 0 \\
			0\arrow{r} & A \arrow{r}{\gamma}\arrow{d} & B' \arrow{r}{\delta}\arrow{d}{\til\psi} & C \arrow{d}{\psi}\arrow{r} & 0 \\
			0\arrow{r} & A \arrow{r}{\alpha} & B \arrow{r}{\beta} & C \arrow{r} & 0  \\
		\end{tikzcd}
	\end{eq}
where $\varphi$ and $\psi$ are isomorphisms. Then there exists an isomorphism $\til\varphi'\colon B \cong B'$ which lifts $\varphi$.
\end{lemm}
With this lemma in combination with \ref{corr}, we have proven
\begin{theo}\label{main}
	Let $\Lambda$ be a noetherian local ring. Given projective birational morphisms $f\colon X\to Y$ and $f'\colon X'\to Y$ of complex algebraic varieties. Suppose that $f$ and $f'$ are semismall and that $X$ and $X'$ are smooth varieties K-equivalent through ${f'}^{-1}\circ f$ (e.g. when $X\to Y \leftarrow X'$ is a semismall flop). Then
	\begin{eq}
		Rf_* \Lambda[n] \cong Rf'_* \Lambda[n].
	\end{eq}
	Suppose further that $X$ and $X'$ are projective varieties, then there is an isomorphism of motives in coefficient $\Lambda$
	\begin{eq}
		M(X)_\Lambda\cong M(X')_\Lambda.
	\end{eq}
\end{theo}
\begin{coro}\label{tor}
	Under the assumptions of \ref{main}, there are isomorphisms of singular cohomology groups
	\begin{eq}
		H^*(X, \bfZ) \cong H^*(X', \bfZ)
	\end{eq}
	and
	\begin{eq}
		H^*(X, \bfF_p) \cong H^*(X', \bfF_p)
	\end{eq}
	for all prime $p$.
\end{coro}
\begin{proof}
	Since $H^k(X, \bfZ)$ and $H^k(X', \bfZ)$ are finitely generated abelian groups, it suffices to compare the ranks and the $p$-primary components $H^k(X, \bfZ)\left[ p^\infty \right]$ and $H^k(X', \bfZ)\left[ p^\infty \right]$ for all primes $p$. \ref{main} provides in particular $H^k(X, \Lambda) \cong H^k(X', \Lambda)$ for all local rings $\Lambda$. Taking $\Lambda = \bfQ$, we see that $\rank_{\bfZ} H^k(X, \bfZ) = \dim_{\bfQ} H^k(X, \bfQ) = \dim_{\bfQ} H^k(X', \bfQ) = \rank_{\bfZ} H^k(X', \bfZ)$. Taking $\Lambda = \bfZ_p$ the ring of $p$-adic integers, we have on the other hand $H^k(X, \bfZ)\left[ p^{\infty} \right]\cong H^k(X, \bfZ_p)\left[ p^{\infty} \right]\cong H^k(X', \bfZ_p)\left[ p^{\infty} \right] \cong H^k(X, \bfZ)\left[ p^{\infty} \right]$. This proves $H^k(X, \bfZ) \cong H^k(X', \bfZ)$. The latter statement is similarly proven by taking $\Lambda = \bfF_p$.
\end{proof}

\section{Proof of \ref{cat}}
In this section we work abstractly on a noetherian abelian category $\calC$ \emph{locally finite} over a noetherian local ring $\Lambda$. This means that the noetherian abelian category $\calC$ has on every Hom set $\Hom_{\calC}(A, B)$ a finitely generated $\Lambda$-module structure for every $A, B\in \calC$, such that each composition map 
\begin{eq}
	\Hom_{\calC}(B, C) \times \Hom_{\calC}(A, B) \to \Hom_{\calC}(A, C)
\end{eq}
is $\Lambda$-bilinear for every $A, B, C\in \calC$.
\smallskip

Given objects $A, C\in \calC$ such that $\Hom_{\calC}(A, C) = 0$, we shall study the $\Lambda$-module $\Ext^1_{\calC}(C, A)$ of extensions of $C$ by $A$. Recall that elements of $\Ext^1_{\calC}(C, A)$ are short exact sequences in $\calC$
\begin{eq}
	0 \to  A \xrightarrow{\alpha} B \xrightarrow{\beta} C \to 0,
\end{eq}
written as $(B, \alpha, \beta)$, modulo the following equivalence relation: we say $(B, \alpha, \beta) \sim (B', \gamma, \delta)$ if there exists a commutative diagramme
\begin{eq}
	\begin{tikzcd}
		0 \arrow{r} & A \arrow{r}{\alpha}\arrow[equals]{d} & B\arrow{r}{\beta}\arrow{d} & C\arrow{r}\arrow[equals]{d} & 0 \\
		0 \arrow{r} & A \arrow{r}{\gamma} & B'\arrow{r}{\delta} & C\arrow{r} & 0
	\end{tikzcd}
	.
\end{eq}
Let $R = \End_{\calC}(A)$ and $S = \End_{\calC}(C)$. Then $\Ext^1_{\calC}(C, A)$ is endowed with a natural $R-S$-bimodule structure in the following manner: For any $\varphi\in R$, the element $\varphi\cdot(B, \alpha, \beta)$ is defined as the second row of the following diagramme
\begin{eq}
	\begin{tikzcd}
		0 \arrow{r} & A \arrow{r}{\alpha}\arrow{d}{\varphi} & B\arrow{r}{\beta}\arrow{d}{(0, 1)} & C\arrow{r}\arrow[equals]{d} & 0 \\
		0 \arrow{r} & A \arrow{r}{(1, 0)} & A\sqcup^{(\varphi, \alpha)}B \arrow{r}{(0,\beta)} & C\arrow{r} & 0
	\end{tikzcd}
\end{eq}
(therein the left square is cocartesian), whereas for $\psi \in S$, the element $(B, \alpha, \beta)\cdot \psi$ is defined as the first row of the following diagramme
\begin{eq}
	\begin{tikzcd}
		0 \arrow{r} & A \arrow{r}{(\alpha, 0)}\arrow[equals]{d} & B\times_{(\beta, \psi)}C\arrow{r}{(0, 1)}\arrow{d}{(1, 0)} & C\arrow{r}\arrow{d}{\psi} & 0 \\
		0 \arrow{r} & A \arrow{r}{\alpha} & B \arrow{r}{\beta} & C\arrow{r} & 0
	\end{tikzcd}
\end{eq}
(therein the right square is cartesian). Here is simple facts concerning such diagrammes
\begin{lemm}
	\begin{enumerate}
		\item[(a)]
			Suppose we have a diagramme of short exact sequences
			\begin{eq}
				\begin{tikzcd}
					0 \arrow{r} & A \arrow{r}{\alpha}\arrow{d}{\varphi} & B\arrow{r}{\beta} \arrow{d}{\chi} & C\arrow{r}\arrow[equals]{d} & 0 \\
					0 \arrow{r} & A' \arrow{r}{\gamma} & B' \arrow{r}{\delta} & C\arrow{r} & 0
				\end{tikzcd}
			\end{eq}
			then the left square is cocartesian.
		\item[(a bis)]
			Dually, suppose we have a diagramme of short exact sequences
			\begin{eq}
				\begin{tikzcd}
					0 \arrow{r} & A \arrow{r}{\alpha}\arrow[equals]{d} & B\arrow{r}{\beta} \arrow{d}{\chi} & C\arrow{r}\arrow{d}{\psi} & 0 \\
					0 \arrow{r} & A \arrow{r}{\gamma} & B' \arrow{r}{\delta} & C'\arrow{r} & 0
				\end{tikzcd}
			\end{eq}
			then the right square is cartesian.
		\item[(b)]
			Given any diagramme of short exact sequence
			\begin{eq}
				\begin{tikzcd}
					0 \arrow{r} & A \arrow{r}{\alpha}\arrow{d}{\varphi} & B\arrow{r}{\beta} \arrow{d}{\chi} & C\arrow{r}\arrow{d}{\psi} & 0 \\
					0 \arrow{r} & A' \arrow{r}{\gamma} & B' \arrow{r}{\delta} & C'\arrow{r} & 0
				\end{tikzcd}
			\end{eq}
			there is a 3-step factorisation
			\begin{eq}
				\begin{tikzcd}
					0 \arrow{r} & A \arrow{r}{\alpha}\arrow{d}{\varphi} & B\arrow{r}{\beta} \arrow{d}{(0, 1)} & C\arrow{r}\arrow[equals]{d} & 0 \\
					0 \arrow{r} & A' \arrow{r}{(1, 0)}\arrow[equals]{d} & A' \sqcup^{(\varphi, \alpha)}B \arrow{r}{(0, \beta)}\arrow{d}{\lambda} & C\arrow{r}\arrow[equals]{d} & 0\\
					0 \arrow{r} & A' \arrow{r}{(\gamma, 0)}\arrow[equals]{d} & B'\times_{(\delta, \psi)} C\arrow{r}{(0, 1)} \arrow{d}{(1, 0)} & C\arrow{r}\arrow{d}{\psi} & 0 \\
					0 \arrow{r} & A' \arrow{r}{\gamma} & B' \arrow{r}{\delta} & C'\arrow{r} & 0
				\end{tikzcd}
			\end{eq}
			where $\lambda\colon A' \sqcup^{(\varphi, \alpha)}B  \to B'\times_{(\delta, \psi)} C$ is defined by 
			\begin{eq}
				\lambda(a', b) = (\gamma(a') + \chi(b), \beta(b)).
			\end{eq}
	\end{enumerate}
\end{lemm}
\begin{proof}
	For (a), consider the following diagramme
	\begin{eq}
		\begin{tikzcd}
			& 0\arrow{d} & 0\arrow{d}& 0\arrow{d} &  \\
			0\arrow{r} & 0 \arrow{r} \arrow{d} & A' \arrow[equals]{r}\arrow{d}{(1, 0)} & A' \arrow{d}{\gamma} \arrow{r}  &0 \\
			0\arrow{r} &  A \arrow{r}{(\varphi, -\alpha)}\arrow[equals]{d} & A' \oplus B\arrow{r}{(\gamma, \chi)}\arrow{d}{(0, 1)} & B'\arrow{d}{\delta}\arrow{r}  & 0 \\
			0\arrow{r} &  A \arrow{r}{\alpha}\arrow{d} & B \arrow{r}{\beta}\arrow{d} & C \arrow{r}\arrow{d} & 0 \\
			& 0 & 0 & 0 &
		\end{tikzcd}
	\end{eq}
	The three columns and the first and the third rows are exact. By 9-lemma, the middle row is also exact. This proves that $(\alpha, \varphi, \chi, \gamma)$ forms a cocartesian square. (a bis) is proven similarly.
	\smallskip

	For (b), it is trivial to check that $\lambda$ is a well defined morphism.
\end{proof}
In our case, $A = A'$ and $C = C'$, so the lemma says that given any commutative diagramme
\begin{eq}
	\begin{tikzcd}
		0 \arrow{r} & A \arrow{r}{\alpha}\arrow{d}{\varphi} & B\arrow{r}{\beta} \arrow{d}{\chi} & C\arrow{r}\arrow{d}{\psi} & 0 \\
		0 \arrow{r} & A \arrow{r}{\gamma} & B' \arrow{r}{\delta} & C\arrow{r} & 0
	\end{tikzcd}
\end{eq}
we always have $\varphi\cdot (B, \alpha, \beta) \sim (B', \gamma, \delta) \cdot \psi$, or $\varphi\cdot [(B, \alpha, \beta)] = [(B', \gamma, \delta)] \cdot \psi$ in $\Ext^1_{\calC}(C, A)$. Of parcitular insterest is the case of \ref{cat}, where $\psi$ is an isomorphism.
\smallskip

Under the setting of \ref{cat}, we have
\begin{eq}
	\begin{tikzcd}
		0\arrow{r} & A \arrow{r}{\alpha}\arrow{d}{\ha \varphi} & B \arrow{r}{\beta}\arrow{d}{\til\varphi} & C\arrow{d}{\varphi}\arrow{r} & 0 \\
		0\arrow{r} & A \arrow{r}{\gamma}\arrow{d}{\ha \psi} & B' \arrow{r}{\delta}\arrow{d}{\til\psi} & C \arrow{d}{\psi}\arrow{r} & 0 \\
		0\arrow{r} & A \arrow{r}{\alpha} & B \arrow{r}{\beta} & C \arrow{r} & 0  \\
	\end{tikzcd}
\end{eq}
with $\psi = \varphi^{-1}$. In other words, this amounts to $\ha\varphi\cdot [(B, \alpha, \beta)] = [(B', \gamma, \delta)]\cdot \varphi$ and $\ha\psi\cdot [(B', \gamma, \delta)] = [(B, \alpha, \beta)]\cdot \varphi^{-1}$ in $\Ext^1_{\calC}(C, A)$. We will invoke the following easy lemma.
\begin{lemm}\label{lift}
	Let $\Lambda$ be a noetherian commutative local ring, and $A$ be an (associative unital) $\Lambda$-algebra which is a finitely generate $\Lambda$-module. Given any surjection of $\Lambda$-algebras $\rho\colon A \to B$, the induced group homomorphism $A^{\times} \to B^{\times}$ is surjective as well.
\end{lemm}
\begin{proof}
	We prove firstly that $\frakm A\subseteq J(A)$, where $J(A)$ is the Jacobson radical of $A$, and similarly that $\frakm B\subseteq J(B)$. The Jacobson radical is the intersection of all maximal left ideals of $A$. Thus it suffices to show that for every maximal left ideal $Q\subseteq A$, the pullback $P$ in $\Lambda$ is equal to $\frakm$, or equivalently that $\Lambda / P$ is a field. To this effect, consider the irreducible left $A$-module $A / Q$, whose endomorphism algebra $D = \End_{A}(A / Q)$ is a finitely generated division $(\Lambda / P)$-algebra. We identify $\Lambda / P$ with a subring of the centre of $D$. Let $0\neq x\in \Lambda / P$ be an element. We shall find its inverse element in $\Lambda / P$. Firstly, $x$ has a two-sided inverse $x^{-1}$ in $D$. Since $D$ is a finite $(\Lambda / P)$-algebra, every element is integral over $\Lambda / P$. Let $(x^{-1})^n + a_1 (x^{-1})^{n-1} + \ldots + a_n = 0$ be a monic polynomial relation for $x^{-1}$, with coefficients $a_i\in \Lambda / P$. Multiplying it with $x^{n-1}$, we obtain $x^{-1} = -\left( a_1 + \ldots + a_n x^{n-1} \right)\in \Lambda / P$. Hence $x^{-1}$ is a inverse of $x$ in $\Lambda / P$, so $\Lambda / P$ is a field. The inclusion $\frakm A\subseteq J(A)$ shows in particular that $A / J(A)$ is a quotient of $A / \frakm A$, thus a semisimple artinian $(\Lambda / \frakm)$-algebra.
	\smallskip

	Now given $b\in B^\times$, the reduction $\ba b\in B / J(B)$ is invertible as well. Since $A / J(A) \to B / J(B)$ is an epimorphism of semisimple artinian algebras, there clearly exists $a\in A$ such that $\ba a\in (A/J(A))^\times$ and $\ba{\rho(a)} = \ba b$. Since $\rho^{-1}(J(B)) = J(A)$, there exists $j\in J(A)$ such that $\rho(a + j) = b$. It remains to verify that $a + j\in A^\times$. Let $a'\in A$ such that $\ba a \ba{a'} = \ba{a'}\ba{a} = \ba{1}\in A / J(A)$. Then $(a + j)a' \in 1 + J(A)$, and that $(a + j)a'$ is invertible according to the characterisation. Similarly, $a'(a + j)$ is invertible. Therefore $a + j\in A^\times$.
\end{proof}

Let $M = R\cdot [(B, \alpha, \beta)]$ be the sub-$R$-module generated by $[(B, \alpha, \beta)]$. Then $M = R\cdot[(B', \gamma, \delta)]\cdot \varphi$. Thus $\ha\varphi$ induces an automorphism on $M$. Put $A = R$ and $B = \image(A\to \End_{\Lambda}(M))$. Then the lemma gives an invertible $\ha\varphi'\in A^\times = \Aut_{\calC}(A)$ such that $\ha\varphi'\cdot[(B, \alpha, \beta)] = [(B', \gamma, \delta)]$. Therefore there is a diagramme of exact sequences
\begin{eq}
	\begin{tikzcd}
		0\arrow{r} & A \arrow{r}{\alpha}\arrow{d}{\ha \varphi'} & B \arrow{r}{\beta}\arrow{d}{\til\varphi'} & C\arrow{d}{\varphi}\arrow{r} & 0 \\
		0\arrow{r} & A \arrow{r}{\gamma} & B' \arrow{r}{\delta} & C \arrow{r} & 0 \\
	\end{tikzcd}
\end{eq}
Since $\varphi$ and $\ha\varphi'$ are both isomorphism, the morphism $B\to B'$ in the middle is also an isomorphism. Thereby have we completed the proof of \ref{cat}.


\end{document}